\documentclass[12pt,emtex]{article}
\usepackage{times}

\usepackage{graphicx}
\usepackage{amsfonts}       %
\usepackage{amsmath} %[leqno]=Gleichungen werden links numeriert;
\usepackage{amssymb}        %
\usepackage{amstext}        %
\usepackage{amsthm}         %theorem-umgebungen
\usepackage{epsfig}
\psfigdriver{dvips}

\usepackage{latexsym}
\usepackage[matrix,curve]{xypic}

\usepackage{rotating}
\rotdriver{dvips}
\usepackage[german]{babel}  %Deutsch
\newcommand{\A}{\mathbb{A}}

\newcommand{\Z}{\mathbb{Z}}
\newcommand{\Q}{\mathbb{Q}}
\newcommand{\QL}{\overline{\mathbb{Q}}_l}
\newcommand{\R}{\mathbb{R}}
\newcommand{\C}{\mathbb{C}}

\begin{document}

\centerline{\bf \Large Existence of Whittaker models}

\bigskip\noindent
\centerline{\bf \Large related to }

\bigskip\noindent
\centerline{\bf \Large four dimensional symplectic Galois
representations}

\bigskip\noindent
\centerline{Rainer Weissauer}

%\newpage

\bigskip\noindent

\bigskip\noindent
%\samepage

\bigskip\noindent Let $\A=\A_{fin}\otimes \R$ be the ring of rational adeles and $GSp(4)$ be the group of symplectic similitudes in four variables. Suppose
$\Pi\cong \Pi_{fin}\otimes\Pi_\infty$ is a cuspidal irreducible automorphic representation of
the group $GSp(4,\A)$, where $\Pi_\infty$ belongs to a discrete series representation of the
group $GSp(4,\R)$. The discrete series representations of the group $GSp(4,\R)$ are grouped
into local $L$-packets \cite{Lec}\cite{Ast}, which have cardinality two and consist of the
class of a holomorphic and the class of a nonholomorphic discrete series representation. Two
irreducible automorphic representations $\Pi=\otimes_v \Pi_v$ and $\Pi'=\otimes_v \Pi'_v$ are
said to be weakly equivalent, if $\Pi_v\cong \Pi'_v$ holds for almost all places $v$.

\bigskip\noindent
The aim of this article is to prove the following

\bigskip\noindent
{\bf Theorem 1}. {\it Let $\Pi$ be a cuspidal irreducible automorphic representation of the
group $GSp(4,\A)$. Suppose $\Pi$ is not CAP and suppose $\Pi_{\infty}$ belongs to the discrete
series of the group $GSp(4,\R)$. Then $\Pi$ is weakly equivalent to an irreducible globally
generic cuspidal automorphic representation $\Pi_{gen}$ of the group $GSp(4,\A)$, whose
archimedean component $\Pi_{gen,\infty}$ is the nonholomorphic discrete series representation
contained in the local archimedean $L$-packet of $\Pi_{\infty}$.}

\bigskip\noindent
For an automorphic representation $\Pi$ as in theorem 1 by \cite{Ast} there exists an
associated Galois representation $\rho_{\Pi,\lambda}$.  From \cite{Ast} and theorem 1 we
obtain, that these Galois representations $\rho_{\Pi,\lambda}$ are symplectic in the following
sense

\bigskip\noindent
{\bf Theorem 2}. {\it Suppose $\Pi$ is as in theorem 1. Then the associated Galois
representation $\rho_{\Pi,\lambda}$ preserves a nondegenerate symplectic $\QL$-bilinear form
$\langle . ,. \rangle$, such that the Galois group acts with the multiplier
$\omega_{\Pi}\mu_l^{-w}$
$$ \langle\ \rho_{\Pi,\lambda}(g) v , \rho_{\Pi,\lambda}(g) w\
\rangle \ = \ \omega_{\Pi}(g)\mu_l^{- w}(g)\cdot \ \langle v, w\rangle \quad , \quad g \in
Gal(\overline{\Q}/\Q) \
$$ where $\mu_l$ is the cyclotomic character.}

\bigskip\noindent
A variant of this construction also yields certain orthogonal four dimensional Galois
representations. See the remark at the end of this article.

\bigskip\noindent
For  an irreducible cuspidal automorphic representation $\Pi$ of $GSp(4,\A)$, which is not CAP
and whose archimedean component belongs to the discrete series, we want to show that $\Pi$ is
weakly equivalent to a globally generic representation $\Pi_{gen}$, whose archimedean component
again belongs to the discrete series. $\Pi$ not being CAP implies, that $\Pi_{gen}$ again is
cuspidal. Hence \cite{Ast}, theorem III  now also holds unconditionally, since the multiplicity
one theorem is known for the generic representation $\Pi_{gen}$.

\bigskip\noindent A careful analysis of the proof shows, that the arguments
imply more. For this we refer to the forthcoming work of U.Weselmann \cite{Wes}.

\bigskip\noindent
\underbar{Proof of theorem 1}: The proof will be based on the hypotheses A,B of \cite{Ast}
proved in \cite{Lec}, and theorem 3 and theorem 4, which will be formulated further below
during the proof of theorem 1. Theorem 3 is a consequence of results of \cite{GRS}. Theorem 4
is proved in \cite{Wes}. For the following it is important, that under the assumptions made in
theorem 1 Ramanujan's conjecture holds for the representation $\Pi$ at almost all places, as
explained in \cite{Ast} section 1.

\bigskip\noindent
{\it Restriction to $Sp(4)$}. The restriction of $\Pi$ to $Sp(4,\A)$ contains an irreducible
constituent, say $\tilde\Pi$. In the notation of \cite{Ast}, section 3 consider the degree five
standard $L$-series
$$\zeta^S(\Pi,\chi,s)\ $$ of $\tilde\Pi$ for $Sp(4)$. For our purposes it suffices to consider this $L$-series
for the primes outside a sufficiently large finite set of place $S$ containing all ramified
places. So this partial $L$-series depends only on $\Pi$, and does not depend on the chosen
$\tilde\Pi$.

\bigskip\noindent
{\it Euler characteristics}.  $\Pi$ is cohomological in the sense of \cite{Ast},\cite{Lec},
i.e. $\Pi$ occurs in the cohomology of the Shimura variety $M$ of principally polarized abelian
varieties of genus two for a suitable chosen $\overline \Q_l$-adic coefficient system ${\cal
V}_{\mu}(\QL)$. Since $\Pi$ is of cohomological type and since we excluded CAP-representations
$\Pi$, the representations $\Pi_{fin}$ and also $\Pi^S$ (i.e. $\Pi^S =\otimes_{v\notin S}
\Pi_v$ outside a finite set $S$ of bad primes) only contribute to the cohomology $H^i(M,{\cal
V}_{\mu}(\QL))$ for the coefficient system ${\cal V}_{\mu}$ \cite{Ast} in the middle degree
$i=3$ and not for the other degrees (see \cite{Ast} hypothesis B(1),(2) and \cite{Lec}). This
cohomological property is inherited to the subgroup $Sp(4,\A)$ by restriction and induction
using the following easy observation: Given two irreducible automorphic representations
$\Pi_1,\Pi_2$ of $GSp(4,\A)$ having a common irreducible constituent after restriction to
$Sp(4,\A)$. Then, if $\Pi_1$ is cuspidal but not CAP, then also $\Pi_2$ is cuspidal and not
CAP.

\bigskip\noindent
Therefore, if we consider the generalized $\Pi^S$-isotypic subspaces for either of the groups
$GSp(4,\A^S)$ or $Sp(4,\A^S)$) in the middle cohomology group of degree $3$, we may as well
replace the middle cohomology group by the virtual representation $H^{\bullet}(M,{\cal
V}_{\mu}(\QL))$. Up to a minus sign this does not change the traces of Hecke operators , which
are later considered in theorem 4. Here $S$ may be any finite set containing the archimedean
place.
%The trace of
%$GSp(4,\A_{fin})$ on $H^{\bullet}(M,{\cal V}_{\mu}(\QL))$ on the
%other hand is rather well behaved.
\goodbreak

\begin{samepage}
\bigskip\noindent
{\it First temporary assumption}. For the moment suppose $\tilde \Pi$ admits a weak lift
$\tilde \Pi'$ to an irreducible automorphic representation of the group $PGl(5,\A)$ in the
sense below (we later show using theorem 4 that this in fact always holds). A representation
$\tilde\Pi'$ of $PGl(5,\A)$ can be considered to be a representation of $Gl(5,\A)$ with trivial
central character. In this sense the lifting property just means, that there exists an
irreducible automorphic representation $\tilde\Pi'$ of $PGl(5,\A)$, for which
$$L^S(\tilde \Pi'\otimes\chi,s) = \zeta^S(\Pi,\chi,s)$$ holds for the standard $L$-series of
$Gl(5)$ and all idele class characters $\chi$ and certain sufficiently large finite sets
$S=S(\chi,\tilde\Pi')$ of exceptional places.

\bigskip\noindent
Considered as an irreducible automorphic representation of $Gl(5,\A)$ the representation
$\tilde\Pi'$ need not be cuspidal (by the way this is  essentially used in the additional
remark on orthogonal representations made at the end after the proof of theorem 1). There exist
irreducible cuspidal automorphic representations $\sigma_i$ of groups $Gl(n_i,\A)$ where
$\sum_i n_i=5$, such that $\tilde\Pi'$ is a constituent of the representation induced from a
parabolic subgroup with Levi subgroup $\prod_i Gl(n_i,\A)$. See \cite{L}, prop.2. Each of the
cuspidal representations $\sigma_i$ can be written in the form $\sigma_i = \chi_i\otimes
\sigma_i^0$ for some unitary cuspidal representations $\sigma_i^0$ and certain one dimensional
characters $\chi_i$.

\bigskip\noindent
Let $\omega_{\sigma_i}$ denote the central character of $\sigma_i$. The identity $L^S(\tilde
\Pi',s) = \zeta^S(\Pi,1,s)$ implies, that the characters
$\omega_{\sigma_i}=\chi^{n_i}\omega_{\sigma_i^0}$ are unitary. In fact, since $\Pi^S$ satisfies
the Ramanujan conjecture by \cite{Ast}, they have absolute value one at all places outside $S$.
Therefore by the approximation theorem all $\chi_i$ are unitary. Hence the $\sigma_i$ itself
have been cuspidal unitary representations.

\bigskip\noindent
Now, since the $\sigma_i$ are cuspidal unitary, the well known theorems of Jaquet-Shalika and
Shahidi on $L$-series for the general linear group \cite{JS} imply the non-vanishing
$L^S(\tilde\Pi'\otimes\chi,1)= \prod_i L^S(\sigma_i\otimes\chi,1)\neq 0$ for arbitrary unitary
idele class characters $\chi$. By the (temporary) assumption, that $\tilde\Pi'$ is a lift of
$\tilde\Pi$, this implies
$$ \zeta^S(\Pi,\chi,1) \ \neq \ 0 \ $$ for all unitary idele class
characters $\chi$.

\bigskip\noindent
{\it Second temporary assumption}. Now in addition we suppose, that $\Pi$ can be weakly lifted
to an irreducible automorphic representation $(\Pi',\omega)$ of the group $Gl(4,\A)\times\A^*$.
By this we mean, that there exists an irreducible automorphic representation $\Pi'$ of
$Gl(4,\A)$ and an idele class character $\omega$, such that for the central characters of $\Pi$
and $\Pi'$ we have
$$\omega_{\Pi'}=\omega^2 \quad \hbox{ and } \quad
\omega_{\Pi}=\omega \ ,$$ and that furthermore
$$L^S(\Pi\otimes\chi,s) = L^S(\Pi'\otimes\chi,s)$$ holds for
sufficiently large finite sets of places $S$ containing all ramified places. Here, following
the notation of \cite{Ast} section three, $L^S(\Pi,s)$ denotes the standard degree four
$L$-series of $\Pi$.
\end{samepage}

\bigskip\noindent
These conditions imposed at almost all unramified places of course completely determines the
automorphic representation $\Pi'$ by the strong multiplicity one theorem for $Gl(n)$. In
particular this implies, that the global lift $\Pi \mapsto (\Pi',\omega)$ commutes with
character twists $\Pi\otimes\chi \mapsto (\Pi'\otimes\chi,\omega\chi^2)$. Furthermore it
implies $(\Pi')^{\vee} \cong \Pi' \otimes \omega^{-1}$. In particular this finally holds
locally at all places including the archimedean place.
%Notice $\Pi_{\infty}\otimes sign
%\cong Pi_{\infty}$ holds for the archimedean discrete series
%representations of $GSp(4,\R)$.

\bigskip\noindent
{\it The archimedean place}.  From the last observation we obtain $\Pi'_\infty\otimes sign
\cong \Pi'_\infty$ from the corresponding $\Pi_\infty\otimes sign \cong \Pi_\infty$, which is
known for discrete series representations $\Pi_\infty$ of $GSp(4,\R)$. Hence it is easy to see,
that $\Pi'_\infty$ respectively $\Pi_\infty$ are determined by their central characters and
their restriction to $Sl(4,\R)$ respectively $Sp(4,\R)$. In fact we will later show, that the
lift $\Pi'$ can be assumed to have a certain explicitly prescribed behavior at the archimedean
place. What this means will become clear later.

\bigskip\noindent
{\it Properties of $\Pi'$}. Again the irreducible automorphic representation $\Pi'$ of
$Gl(4,\A)$ need not be cuspidal. $\Pi'$ is a constituent of a representation induced from some
parabolic subgroup with Levi subgroup $\prod_j Gl(m_j)$ for $\sum_j m_j=4$, with respect to
some irreducible cuspidal representations $\tau_j$ of $Gl(m_j,\A)$.

\bigskip\noindent
The same argument, as already used for the first temporary assumption, implies the $\tau_j$ to
be unitary cuspidal. This excludes $m_j=1$ for some $j$, since otherwise this would force the
existence of a pole of $L^S(\Pi\otimes \chi,s)$ for $\chi=\tau_j^{-1}$ at $s=1$ again by the
results \cite{JS} of Jaquet-Shalika and Shahidi on the analytic behavior of $L$-series for
$Gl(n)$ on the line $Re(s)=1$. Notice, a pole at $s=1$ would imply $\Pi$ to be a CAP
representation (see \cite{P}). This however contradicts the assumptions of theorem 1, by which
$\Pi$ is not a CAP representation.

\bigskip\noindent
{\it Let us return to $\Pi$}. Either $\Pi'$ is cuspidal; or $\Pi'$ comes by induction from a
pair $(\tau_1,\tau_2)$ of irreducible cuspidal representations $\tau_i$ of $Gl(2,\A)$, for
which $\omega_{\tau_1}\omega_{\tau_2}=\omega^2$ holds. In this second case we do have
$$L^S(\Pi',s)=L^S(\tau_1,s)L^S(\tau_2,s)\ .$$ By $L^S(\Pi,s)=L^S(\Pi',s)$ therefore $\Pi$ is
a weak endoscopic lift, provided the central characters $\omega_{\tau_1}=\omega_{\tau_2}$
coincide. For weak endoscopic lifts theorem 1 obviously holds. See \cite{Ast}, hypothesis A
part (2) and (6) combined with \cite{Lec}. To complete the discussion of this case we establish
the required identity $\omega_{\tau_1}=\omega_{\tau_2}$. For this we need some further argument
and therefore we make a digression on theta lifts first.

\bigskip\noindent
{\it The theta lift}. $Gl(1)$ acts on $Gl(4)\times Gl(1)$ such that $t$ maps $(h,x)$ to
$(ht^{-1},xt^2)$. By our temporary assumptions made, the central character of $\Pi'$ is
completely determined so that the representation $(\Pi',\omega)$ of $Gl(4)\times Gl(1)$
descends to a representation on the quotient group $G(\A)=(Gl(4,\A)\times Gl(1,\A))/Gl(1,\A)$.
This quotient group $G(\A)$ is isomorphic to the special orthogonal group of similitudes
$GSO(3,3)(\A)$ attached to the split 6 dimensional Grassmann space $\Lambda^2(\Q^4)$ with the
underlying quadratic form given by the cup-product.

\bigskip\noindent
The (generalized) theta correspondence of the pair
$$\bigl(GSp(4),GO(3,3)\bigr)\ $$ preserves central characters. If we apply
the corresponding theta lift to the representation $\Pi$ of
$GSp(4,\A)$, then according to \cite{AG} p.40 the theta lift of
$\Pi$ to $GO(3,3)(\A)$ is nontrivial if and only if $\Pi$ is a
globally generic representation. In this case it is easy to see,
that the lift is globally generic. For the converse we need the
following result announced by Jaquet, Piateskii-Shapiro and
Shalika \cite{JPS}. See also \cite{AG} and \cite{S2}, \cite{PSS}.

\bigskip\noindent
{\bf Theorem 3}. \label{hypothesis C} {\it An irreducible cuspidal automorphic representation
$\Pi'$ of $Gl(4,\A)$ lifts nontrivially to $GSp(4,\A)$ under the generalized theta
correspondence of the pair $(GSp(4),GO(3,3))$ if and only if the alternating square degree six
$L$-series $L(\Pi',\chi,s,\Lambda^2)$ or equivalently some partial $L$-series
$$L^S(\Pi',\chi,s,\Lambda^2)$$ (for a suitably large finite set of places S containing all bad
places) has a pole at $s=1$ for some unitary idele class character $\chi$. In this case the
lift of $\Pi'$ is globally generic, and also $\Pi$ is generic.}

\bigskip\noindent
{\it Remark on the degree six $L$-series}. To apply this it is enough to observe, that under
our second temporary assumptions we have enough control on $L^S(\Pi',\chi,s,\Lambda^2)$ to
apply theorem 3 for the representation $\Pi'$ of $Gl(4,\A)$, attached to the representation
$\Pi$ of $GSp(4,\A)$ subject to our second temporary assumption. Indeed by an elementary
computation the second temporary assumption implies the identity
$$ L^S(\Pi',\chi,s,\Lambda^2) =
L^S(\omega\chi,s)\zeta^S(\Pi,\omega\chi,s) \ ,$$ where $\omega=\omega_{\Pi}$ denotes the
central character of $\Pi$. Now compare $L^S(\Pi',\chi,s,\Lambda^2)$ with the standard
$L$-series of the special orthogonal group $SO(3,3)$, which is used in the proof of theorem 3:

\bigskip\noindent We  have exact sequences
$$ \xymatrix{1 \ar[r] & SO(3,3) \ar[r] & GSO(3,3) \ar[r]^-{\lambda} & Gl(1)
\ar[r] & 1 \cr} \ ,$$
$$ \xymatrix{1 \ar[r] & Gl(1) \ar[r]^-{i} &  Gl(4)\times Gl(1) \ar[r] & GSO(3,3)
\ar[r] & 1 \cr} \ ,$$ where $\lambda$ is the similitude homomorphism and where
$i(t)=(t^{-1}\cdot id, t^2 )$. Hence for the Langlands dual groups, which at the unramified
places describe the restriction of spherical representations of $GSO(3,3)(F)$ to $SO(3,3)(F)$,
we get
$$  \xymatrix{1 \ar[r] & \C^* \ar[r]^-{\hat\lambda} & \widehat{GSO(3,3)} \ar[r] & \widehat{SO(3,3)}
\ar[r] & 1 \cr} \ .$$ Notice $\widehat{GSO(3,3)} \subseteq \hat Gl(4) \times \hat Gl(1)$. The
6-dimensional complex representation of $Gl(4,\C)\times Gl(1,\C)$ on $\Lambda^2(\C^4)$ defined
by
$$(A,t)\cdot X = t^{-1}\cdot \Lambda^2(A)(X)$$ is trivial on the
subgroup $\hat\lambda(\C^*)$, hence defines a 6-dimensional representation of the $L$-group
$\widehat{SO(3,3)}$. The $L$-series of this representation defines the degree 6 standard
$L$-series $L^S(\sigma,s)$ of an irreducible automorphic representation $\sigma$ of the group
$SO(3,3)(\A)$. Apparently, for $\sigma$ spherical outside $S$ in the restriction of
$(\Pi',\chi)$ so that $\omega_{\Pi'}=\chi^2$, we therefore get
$$ L^S(\sigma,s) = L^S(\Pi',\chi^{-1},s,\Lambda^2) \ .$$

\bigskip\noindent
\goodbreak\noindent
\underbar{Proof of theorem 3}: Using the remark on the degree six $L$-series we now can invoke
\cite{GRS}, theorem 3.4 to deduce hypothesis C.  The condition on genericity made in loc. cit.
automatically holds for an cuspidal irreducible automorphic representation $\sigma$ of
$SO(n,n)(\A)$ for $n=3$, since this conditions is true for $Gl(4,\A)$. Therefore by \cite{GRS}
a pole of $L^S(\sigma,s)$ at $s=1$ implies, that $\sigma$ has a nontrivial, cuspidal generic
theta lift to the group $Sp(2n-2)(\A)= Sp(4,\A)$. This in turn easily implies the same for the
extended theta lift from $GSO(3,3)(\A)$ to $GSp(4,\A)$. This proves theorem 3.

\bigskip\noindent
{\it Additional remark}. If in addition the spherical representation $\Pi'^S$ is a constituent
of an induced representation attached to a pair of unitary cuspidal irreducible automorphic
representations $\tau_1,\tau_2$ of $Gl(2,\A)$, then furthermore
$$ L^S(\Pi',\chi,s,\Lambda^2) \ = \
L^S(\tau_1\times(\tau_2\otimes\chi),s)L^S(\omega_{\tau_1}\chi,s) L^S(\omega_{\tau_2}\chi,s) \
.$$ This, as well as $\omega^2=\omega_{\Pi'}=\omega_{\tau_1}\omega_{\tau_2}$, are rather
obvious. But by the wellknown analytic properties of $L$-series [JS] it implies, that the right
side now has poles for $\chi= \omega_{\tau_i}^{-1}$ and $i=1,2$. Notice the $\tau_i$ are
unitary cuspidal.

\bigskip\noindent
This being said we now complete the discussion of the case, where $\Pi'$ is not cuspidal.

\bigskip\noindent
{\it Reduction to the case $\Pi'$ cuspidal}. If $\Pi'$ is not cuspidal, then as already shown
$\Pi'$ is obtained from a pair of irreducible unitary cuspidal representations
$(\tau_1,\tau_2)$ by induction. To cover theorem 1 in this case we already remarked, that it
suffices to show $\omega_{\tau_1}=\omega_{\tau_2}$. If these two characters were different,
then $\chi_i=\omega\omega_{\tau_i}^{-1}\neq 1$ for $i=1$ or $i=2$.
%But since $\chi_1\chi_2=1$, this must hold for both choices.
Therefore $L^S(\Pi',\chi,s,\Lambda^2)$ would have a pole at $s=1$ for $\chi =
\omega_{\tau_i}^{-1}$ by the previous \lq{additional remark}\rq. Hence the identity $
L^S(\Pi',\omega_{\tau_i}^{-1},s,\Lambda^2) = L^S(\chi_i ,s)\zeta^S(\Pi,\chi_i,s)$ would imply
the existence of poles at $s=1$ for the $L$-series $\zeta^S(\Pi,\chi_i,s)$.
%and both choices of $\chi_i$.
Thm.4.2 of \cite{Ast}, section 4 then would imply $(\chi_i)^2=1$. Since $\chi_1\chi_2=1$ holds
by definition of $\Pi'$, therefore $\chi_1=\chi_2$. Thus $\omega_{\tau_1}=\omega_{\tau_2}$. So
we are in the case already considered in \cite{Ast}: $\Pi$ is a weak lift. In this case the
statement of theorem 1 follows from the multiplicity formula for weak endoscopic lifts
\cite{Ast}, hypothesis A (6). Thus we may suppose from now on, that $\Pi'$ is cuspidal.

\bigskip\noindent
%In the other direction, since $\chi=\omega^{-1}$ produces a pole
%of $L^S(\omega\chi,s)$ at $s=1$ this implies the case the
%existence of a pole of
%$L^S(\tau_1\times\tau_2\otimes\omega^{-1},s)$ at $s=1$. Here we
%used that $\zeta^S$ does not vanish at $s=1$, as explained above.
%As a consequence of this pole we get $\tau_1^{\vee}\cong
%\tau_2\otimes\omega^{-1}$ or $\tau_1\cong \tau_2\otimes
%(\omega_1/\omega)$. Put $\chi_i=\omega_{\tau_i}/\omega$. These
%are the quadratic characters, which define the algebras $K_i$.
{\it Applying theorem 3}. Both our temporary assumptions on the existence of the lifts
$\tilde\Pi'$ and $\Pi'$ imply, that from now on we can assume without restriction of
generality, that $\Pi'$ is a unitary cuspidal representation. Since we deduced
$\zeta^S(\Pi,\chi,1)\neq 0$ for all unitary characters $\chi$ from our first temporary
assumption, the crucial identity $ L^S(\Pi',\chi,s,\Lambda^2) =
L^S(\omega\chi,s)\zeta^S(\Pi,\omega\chi,s) $ forces the existence of a pole for
$L^S(\Pi',\omega^{-1},s,\Lambda^2)$ at $s=1$. Therefore we are in a situation where we can
apply theorem 3: Since $\Pi'$ is cuspidal, the pair $(\Pi',\omega)$ defines a cuspidal
irreducible automorphic representation of $GSO(3,3)(\A)$. It nontrivially gives a backward lift
from $GSO(3,3)$ to a globally generic automorphic representation $\Pi_{gen}$ of $GSp(4,\A)$
using theorem 3. Comparing both lifts at the unramified places gives
$$ L^S(\Pi,s) = L^S(\Pi',s) = L^S(\Pi_{gen},s) \ .$$
The first equality holds by assumption. The second equality follows from the behavior of
spherical representations under the Howe correspondence \cite{R}. See also \cite{PSS}, p.416
for this particular case. Hence $\Pi$ and the generic representation $\Pi_{gen}$ are weakly
equivalent.

\bigskip\noindent
In other words, using two temporary assumptions, we now have almost deduced theorem 1. In fact
the generic representation $\Pi_{gen}$, that has been constructed above, is weakly equivalent
to $\Pi$. Hence it is cuspidal, since $\Pi$ is not CAP. However, for the full statement of
theorem 1 one  also needs control over the archimedean component of $\Pi_{gen}$. We postpone
this archimedean considerations for the moment and rather explain first, how to establish the
two temporary assumptions to hold unconditionally. This will be deduced from the topological
trace formula.

\bigskip\noindent
{\it Construction of the weak lifts $\Pi'$ and $\tilde\Pi'$ of $\Pi$}. The existence of these
lifts will follow from a comparison of the twisted topological trace formula of a group
$(G,\sigma)$ with the ordinary topological trace formula for a group $H$ for the pairs
$H=GSp(4,\A)$ and $(G,\sigma)=(GSO(3,3)(\A),\sigma)$ respectively $\tilde H=Sp(4,\A)$ and
$\tilde G=(PGl(5,\A),\tilde\sigma)$. Here $\sigma$ respectively $\tilde\sigma$ denote
automorphisms of order two of $G$ resp. $\tilde G$. In both cases the group of fixed points in
the center under the automorphism $\sigma, \tilde\sigma$ will be a Zariski connected group, a
condition imposed in [W3]. Notice $G$ is isomorphic to the quotient of $Gl(4)\times Gl(1)$
divided by the subgroup $S$ of all zentral elements of the form $(z\cdot id,z^{-2})$, hence is
isomorphic to $Gl(4)/\{\pm 1\}$.  Hence for a local field $F$
$$ GSO(3,3)(F)\cong \Bigl( Gl(4,F)\times F^*\Bigr) /F^* \ , $$
where $t\in F^*$ acts on $(h,x)\in Gl(4,F)\times F^*$ via $(h,x)\mapsto (ht^{-1},xt^2)$. The
group $GSO(3,3)(F)$ can be realized to act on the six dimensional Grassmann space
$\Lambda^2(F^4)$. See \cite{AG}, p.39ff and \cite{Wa}, p.44f.  This identifies the quotient
group $G(F)$  with the special orthogonal group of similitudes $GSO(3,3)(F)$ attached to the
split 6 dimensional quadratic space $\Lambda^2(F^4)$ defined by the cup-product. The similitude
character is $\lambda(h,x)=det(h)x^2$. The automorphism $\sigma$ of $G$ is induced by the map
$(h,x)\mapsto (\omega (h')^{-1} \omega^{-1}, det(h)x)$ for a suitable matrix $\omega\in
Gl(4,F)$ chosen in such a way, that $\sigma$ stabilizes a fixed splitting \cite{BWW}.1.9. Then
$\sigma$ is the identity on the center $Z(G)\cong Gl(1)$ of $G$.

\bigskip\noindent
Let us start with some

\bigskip\noindent
{\it Notations}. Let $G$ be a split reductive $\Q$-group with a fixed splitting
$(B,T,\{x_{\alpha}\})$ over $\Q$ and center $Z(G)$. Let $\sigma$ be a $\Q$-automorphism of $G$,
which stabilizes the splitting, such that the group  $Z(G)^{\sigma}$ of fixed points is
connected for the Zariski topology. Let $K_{\infty}^+$ be the topological connected component
of a maximal compact subgroup of $G(\R)$, similarly let $Z_{\infty}^+$ be the connected
component of the $\R$-valued points of the maximal $\R$-split subtorus of the center of $G$.
Let $X_G=G(\R)/K_{\infty}^+Z_{\infty}^+$ be the associated symmetric domain as in \cite{Wes}
(5.22). Let $V$ be an irreducible finite dimensional complex representation of $G$ with highest
weight $\chi \in X^*(T)$, which is invariant under $\sigma$. It defines a bundle $
V_G=G(\Q)\setminus [G(\A_{fin})\times X_G\times V]$ over $M_G=G(\Q)\setminus [G(\A_{fin})\times
X_G]$. See also \cite{Wes} (3.4). Let ${\cal V}_\chi$ denote the associated sheaf. For the
natural right action of $G(\A_{fin})$ on $M_G$ and $V_G$ the cohomology groups $H^i(M_G,{\cal
V}_{\chi})$ become admissible $G(\A_{fin})$-modules, on which $\sigma$ acts. Let
$H^{\bullet}(M_G,{\cal V}_{\chi})= \sum_i (-1)^i H^i(M_G,{\cal V}_{\chi})$ be the corresponding
virtual modules. Then for $f\in C_c^{\infty}(G(\A_{fin}))$ the traces
$$ T(f,\sigma,G,\chi) = Trace(f \cdot \sigma,
H^{\bullet}(M_G,{\cal V}_{\chi}))$$ are well defined.

\bigskip\noindent
Let $G_1$ be the maximal or stable endoscopic group for $(G,\sigma,1)$ in the sense of
$\sigma$-twisted endoscopy (see \cite{KS}(2.1) and \cite{Wes}, section 5), where we assume that
the character $\omega=\omega_a$ is trivial. The corresponding endoscopic datum $(G_1,{\cal
H},s,\xi)$ has the property ${\cal H}= {}^L G_1$. The dual group $\hat G_1$ is the group of
fixed points of $\hat G$ under $\hat \sigma$ and is again split, defining $\xi: {}^L G_1 \to
{}^L G$. Let $T_1$ be a maximal $\Q$-torus of $G_1$, then we can identify the character group
$X^*(T_1)$ with the fixed group $X^*(T)^{\sigma}$. Hence $\chi$ defines a coefficient system
${\cal V}_{\chi_1}$ on $M_{G_1}$ and we can similarly define $T(f_1,id,G_1,\chi_1)$. Functions
$f=\prod_{v\neq \infty} f_v\in C_c^{\infty}(G(\A_{fin}))$ and $f_1=\prod_{v\neq \infty}
f_{1,v}\in C_c^{\infty}(G_1(\A_{fin}))$ are called matching functions, if each of the local
pairs are matching in the sense of \cite{KS} (5.5.1) up to $z$-extensions, so that in
particular $f_v$ and $f_{v,1}$ are characteristic functions of suitable hyperspecial maximal
compact subgroups $K_v\subseteq G_v=G(\Q_v)$ resp. $K_{v,1}\subseteq G_{1,v}=G_1(\Q_v)$ for
almost all $v\notin S$ ($S$ a suitable finite set of places which may be chosen arbitrarily
large).  By simplicity here we tacitly neglect, that $f_1$ and $f$ have to be chosen as in
\cite{KS} p. 24 or 70 up to an integration over the central group denoted $Z_1(F)$ in loc. cit.
The functions $f_1=\prod_{v} f_{1,v}$ and $f=\prod_v f_v$ are said to be {\it globally
matching} functions, if the analog of formula \cite{KS} (5.5.1) holds for the global stable
orbital integrals defined over the finite adeles all global elements $\delta,\delta_H$. This
slightly weaker global condition suffices for the comparison of trace formulas.

\bigskip\noindent
{\it The homorphism $b_\xi$}. For $v\notin S$ the classes of irreducible $K_{1,v}$-spherical
representations $\Pi_{v,1}$ of $G_1(\Q_v)$ are parameterized by their \lq{Satake parameter}\rq
$\alpha_1\in Y_1(\C)=\hat T_1/W(\hat G_1)$. Since $Y_1(\C)$ can be identified with $(\hat
G_1)_{ss}/int(\hat G_1)$ (see \cite{Bo}, lemma 6.5) and similar for $G$, the endoscopic map
$\xi$ induces an algebraic morphism $(\hat G_1)_{ss}/int(\hat G_1) \to (\hat G)_{ss}/int(\hat
G)$, hence a map $\xi_Y: Y_1(\C) \to Y(\C)$. Since the spherical Hecke algebra of
$(G_{1,v},K_{1,v})$ resp. $(G_v,K_v)$ can be identified with the ring of regular functions
$\C[Y_1]$ resp. $\C[Y]$, we obtain an induced algebra homomorphism $b_{\xi}=\xi_Y^*$ from the
spherical Hecke algebra $C_c^{\infty}(G_v/\!/K_v)$ to spherical Hecke algebra
$C_c^{\infty}(G_{1,v}/\!/K_{1,v})$. The Satake parameter $\alpha=\xi_Y(\alpha_1)$ defines a
$K_v$-spherical representation of $G_{1,v}$, also denoted $\Pi_{v}=r_{\xi}(\Pi_{v,1})$. By
construction of the map it satisfies $\Pi_v^{\sigma}\cong \Pi_v$. Hence the action of $G_v$ on
$\Pi_v$ can uniquely be extended to an action of the semidirect product $G_v\rtimes<\sigma
>$ assuming, that $\sigma$ fixes the spherical vector.

\bigskip\noindent
{\it The topological trace formula}. We give a review of the $\sigma$-twisted topological trace
formula of Weselmann \cite{Wes} theorem 3.21, which generalizes the topological trace formula
of Goresky and MacPherson \cite{GMP}, \cite{GMP2} from the untwisted to the twisted case.
Weselmann shows that these trace formulas themselves are \lq{stable}\rq\ trace formulas
\cite{Wes} theorem 4.8 and remark 4.5, i.e. can be written entirely in terms of stable twisted
orbital integrals $SO_\gamma^{G,\sigma}(f)$ for the maximal elliptic endoscopic group of
$(G,\sigma)$ in the sense of twisted endoscopy
$$ T(f,\sigma,G,\chi) = \sum_{I\subset \Delta, \sigma(I)=I} (-1)^{\vert (\Delta \setminus
I)/\sigma\vert} \cdot T_I(f,\sigma,G,\chi) $$
$$ T_I(f,\sigma,G,\chi) \ =\ \sum_{\gamma \in P_I(\Q)/\sim}^{'} \alpha_\infty(\gamma,1) \cdot
SO_\gamma^{G,\sigma}(f) \cdot Trace(\gamma\cdot\sigma,V) $$ where the summation runs over all
stable conjugacy classes $\gamma$ in $P_I(\Q)$ of $I$-contractive elements, whose norm is in
$L_\infty^I$ (see loc. cit theorem 4.8). In principle these trace formulas can be compared
without simplifying assumptions (see \cite{Wes}, section 5) except for a different notion of
matching functions due to the factors $\alpha_\infty(\gamma,1)$. See \cite{Wes} (5.15).

\bigskip\noindent
{\it Strongly matching functions}. Globally matching functions $(f,f_1)$  are said to be {\it
strongly matching}, if there exists a universal constant $c=c(G,\sigma)\neq 0$ such that
%there exists a constant $c\neq 0$ such that $
%SO_{\gamma_0}^{G,\sigma}(f) = SO_{\gamma_1}^{G_1}(f_1)=0$ holds for all semisimple strongly
%$\sigma$-regular elements $\gamma_0\in G(\Q)$ respectively matching elements $\gamma_1\in
%G_1(\Q)$ for which
$ T_I(f,\sigma,G,\chi) = c\cdot T_I(f_1,id,G_1,\chi)$ holds for all $I$. Then
$$T(f,\sigma,G,\chi)=  c\cdot T(f_1,id,G_1,\chi)$$ holds by definition. In the
lemma below we will show for constants $c_I=c(G,\sigma,I)\neq 0$ $$\vert
\alpha_\infty(\gamma_0,1)\vert =c_I\cdot \vert \alpha_\infty(\gamma_1,1)\vert$$  for all
summands of the sum defining $T_I(f,\sigma,G,\chi)$ respectively $T_I(f_1,id,G_1,\chi)$ and
sufficiently regular $\gamma_0$ respectively $\gamma_1$.
 Notice
$(\gamma_0,\gamma_1)= (\delta,\delta_H)$ are the global $\Q$-rational elements to be compared
in the notions above. For $\gamma=\gamma_0$ (or $\gamma=\gamma_1$) and the respective group $G$
(or $G_1$) by \cite{Wes} theorem 4.8
$$ \vert\alpha_\infty(\gamma,1)\vert = \frac{O_\sigma^\infty(\gamma,1)\cdot \#
H^1(\R,T)}{d^I_{\zeta,\gamma} \cdot vol_{db_\infty'}(\overline
{(G^I_{\gamma,\sigma})}'(\R)/\tilde\zeta)} \ .$$ This is the factor for the terms in the sum
defining $T_I(f,\sigma,G,\chi)$, i.e. $\gamma\in P_I(\Q)$ for the $\sigma$-stable $\Q$-rational
parabolic subgroup $P_I=M_I U_I$. We now discuss the ingredients of the formula defining
$\alpha_\infty(\gamma,1)$.

\bigskip\noindent
{\it Regularity}. We will later choose a pair of matching functions $(f,f_1)$, whose stable
orbital integrals have strongly $\sigma$-regular semisimple support. For the notion of strongly
$\sigma$-regular see \cite{KS} p.28. Taking this for granted, we therefore analyze the above
terms under the assumption, that $Int(\gamma)\circ \sigma$ is strongly $\sigma$-regular
semisimple. For simplicity of notation we consider the case $(G,\theta)$ in the following. The
case of its maximal elliptic $\sigma$-endoscopic group $(G_1,id)$ of course is analogous. We
will also assume $Z(G)^\theta$ to be Zariski connected, since this suffices for our
applications. Although the following discussion holds more generally, at the end we restrict to
our two relevant cases for convenience.

%We will also choose the matching pair
%of functions, such that only elements $\gamma=\gamma_0,\gamma_1$ need to be considered in the
%trace trace formulas \cite{Wes} theorem 4.8 with a $\Q$-anisotropic twisted centralizer modulo
%the center. This will achieved by a suitable choice of $(f,f_1)$ at one local nonarchimedean
%place. By this choices we can assume the Levi group $L$ to be maximal.

\bigskip\noindent
{\it The groups $G^I_{\gamma,\sigma}$}. Fix a a $\sigma$-stable rational parabolic subgroup
$P_I$ of $G$ (or similar $G_1$). For a strongly $\sigma$-regular semisimple elements
$\theta=Int(\gamma)\circ \sigma$ , where $\gamma\in P_I(\Q)$, the twisted centralizer
$G^I_{\gamma,\sigma}= (P_I)^{\theta}$ is abelian. In fact, $G^{\theta}$ is abelian  by
\cite{KS} p.28 and the centralizer $G^{\theta}$ in $G$ is a maximal torus ${\rm T} \subseteq
G$. ${\rm T}$  is $\theta$-stable and $T:=G_{\gamma,\sigma}$ is the group of $\theta $-fixed
points in ${\rm T}$ $$ G^\theta = {\rm T}^\theta \ .$$ $\theta$ is strongly $\sigma$-regular,
hence there exists a pair $({\rm T},B)$ ($B$ a Borel containing ${\rm T}$ defined over the
algebraic closure), which is $\theta$-stable. Since ${\rm T}$ acts transitively on splittings,
there exists a $t$ in ${\rm T}$ (again over the algebraic closure), such that
$\theta^*=int(t)\theta$ respects a fixed splitting $({\rm T},B,\{x_\alpha\})$. Then
$${\rm T}^\theta = {\rm T}^{int(t)\theta}={\rm T}^{\theta^*}\ .$$ By \cite{KS} p.14
$G^{\theta^*}$ is Zariski connected if and only if ${\rm T}^{\theta^*}$ is Zariski connected.
Therefore ${\rm T}^\theta$ is Zariski connected, if $G^{\theta^*}$ is Zariski connected. Now
$G^{\theta^*} = G^1 \cdot Z(G)^{\theta^*}$ by \cite{KS}, p.14, for the Zariski connected
component $G^1=(G^{\theta^*})^0$. By our assumption $Z(G)^{\theta^*} =Z(G)^{\sigma}$ is Zariski
connected. Hence $G^{\theta^*}$ and therefore ${\rm T}^{\theta}$ are both Zariski connected. In
other words ${\rm T}^{\theta}$ is a subtorus of ${\rm T}$. Notice $\sigma^2=1$ implies
$(\theta^*)^2=1$, since an inner automorphism fixing a splitting is trivial. Therefore to
describe $({\rm T},\theta^*)$ over the algebraic closure, we can replace $\theta^*$ by our
original automorphism $\sigma$, which also fixes (some other but conjugate) splitting of $G$.
Over the algebraic closure $(T,\theta^*)$ is isomorphic to a direct product $\prod_i
(T_i,\theta_i^*)$, where the factors are either $(Gl(1),id)$ or $(Gl(1),inv)$ or
$(Gl(1)^2,\theta^*)$, where $\theta^*$ is the flip automorphism of the two factors
$\theta^*(x,y)=(y,x)$. Which types arise does not depend on $I$ as shown below, and can be
directly read of from the way in which $\sigma$ acts on the $\sigma$-stable diagonal reference
torus \cite{BWW}. Two cases are relevant: \cite{BWW} example 1.8 where $G=PGl_{2n+1}$ with
$\sigma(g) =J ({}^t g^{-1}) J^{-1}$ and $J$ as in loc. cit.,  and \cite{BWW} example 1.9 where
$G=Gl(2n) \times Gl(1)$ and $\theta(g,a)=(J ({}^t g^{-1}) J^{-1},det(g)a)$. For $n=2$ these
specialize to the cases considered in theorem 4, except that in example 1.9 one has to divide
by $Gl(1)$ to obtain $GSO(3,3)$ as explained already.

\bigskip\noindent
{\it Absolute independence from $I$}. ${\rm T}$ is the unique maximal torus of $G$, which
contains $G^\theta$ as subgroup and which is $\theta$-stable. Furthermore $G^I_{\gamma,\sigma}
= {\rm T}^\theta \cap P_I$. Since $U_I\cap {\rm T}^\theta$ is trivial, the projection $P_I\to
P_I/U_I\cong M_I$ induces an isomorphism of  $ G_{\gamma,\sigma}^I$ with its image in $M_I$. We
can find a $\theta$-stable maximal torus in $M_I$ containing the image. By dimension reason,
this maximal torus coincides with the image of  $ G_{\gamma,\sigma}^I$. By the uniqueness of
${\rm T}$ this determines ${\rm T}$ and $T={\rm T}^\theta$ within the subgroup $(M_I,\theta)$
of $(G,\theta)$. For strongly regular $\gamma$ this implies $$ G_{\gamma,\sigma}^I \cong {\rm
T}^\theta = {\rm T}^{\theta^*} = T \ ,$$ which is independent from $I$. In our situation
$\sigma=\eta_1$ holds for $\eta_1$ as defined in \cite{Wes} (2.1). By the Gauss-Bonnet formula
\cite{Wes} (3.9) the torus $T$ is $\R$-anisotropic modulo the center of $M_I$. By a global
approximation argument therefore the groups $\tilde \zeta$ in the formula above are trivial.
Similarly, since in our case the centers of $G$ and $G_1$, hence also the centers of their
respective $\Q$-Levi subgroups $M_I$, are split $\Q$-tori, the groups $\zeta$ can be assumed to
be trivial. Hence $d^I_{\zeta,\gamma}=1$ by \cite{Wes} (2.23).

\bigskip\noindent
{\it The other factors}. By \cite{Wes} Lemma (2.17) and (2.15) and the isomorphism
$\pi_0(\tilde F(g_\eta,\gamma)^{\tilde\eta_{\gamma,h_\infty}}) \cong \pi_0(\tilde
F(g_\eta,\gamma))^{\tilde\eta_{\gamma,h_\infty}}$ the factors $O_\sigma^\infty(\gamma,1) = \#
R_{\gamma,\sigma}$ are
$$ O_\sigma^\infty(\gamma,1) = \frac{\#\Bigl( K_\infty^{I,m}/K_\infty^I\Bigr)^{\eta_2}}
{ \# \pi_0\Bigl( G_{\gamma,\sigma}^I(\R) / (G_{\gamma,\sigma}(\R) \cap p_1 K_\infty^+
Z_\infty^+ A_I p_1^{-1}) \Bigr) } \ .$$  $G_{\gamma,\sigma}^I(\R)$ is topologically connected,
since $G_{\gamma,\sigma}^I=T$ is a torus. Therefore the denominator is trivial, and the factor
$ O_\sigma^\infty(\gamma,1)$ becomes $( K_\infty^{I,m}/K_\infty^I)^{\eta_2}$ with notations
from \cite{Wes}, (2.2). It only depends on $I$, but does not depend on $\gamma$. So we need
compare the factors
$$ \vert \alpha_\infty(\gamma,1)\vert \ =\ \#\bigl( K_\infty^{I,m}/K_\infty^I\bigr)^{\eta_2} \cdot \frac{ \# H^1(\R,T) }{ vol_{db_\infty'}(T'(\R))} \ .$$
The $\sigma$-stable torus ${\rm T}$ contains  $G_{\gamma,\sigma}^I=T ={\rm T}^{\sigma}$
($\sigma$-invariant subtorus). Its $\R$-structure might  a priori depend on  $\gamma$, but in
fact does not. $T'=({\rm T}^\sigma)'$ is the maximal $\R$-anisotropic subtorus of $T$, as
follows from the description of \cite{Wes} (3.9).
%and the fact already mentioned above, that $T$ is $\R$-anisotropic modulo
%the center of $M_I$.
Now, since
$$ G_{\gamma_0,\sigma}^I=T\cong {\rm T}^\sigma\  \longrightarrow\  {\rm T}_\sigma \cong T_1=(G_1)_{\gamma_1,id}^{I_1} $$
are isogenious tori, the definition of measures $db_\infty'$ (\cite{Wes}, (3.9)) and the
measures used the definition of matching of Kottwitz-Shelstad (\cite{KS}, p. 71) shows, that
both factors $ vol_{db_\infty'}((G_{\gamma_0,\sigma}^I)'(\R))$ for $G$ and $\gamma_0$ and $
vol_{db_\infty'}((G_1)_{\gamma_1,id}^I)'(\R))$ for $G_1$ and $\gamma_1$ differ by a constant
independent of $(\gamma_0,\gamma_1)$. Recall that $\gamma=\gamma_0\in G(\Q)$ and
$\gamma=\gamma_1\in G_1(\Q)$ is a pair of points related by the Kottwitz-Shelstad norm. For
such a pair we have the tori $T\subseteq G$ and $T_1 \subseteq G_1$ and $T_1 = {{\rm
T}_\sigma}$ ($\sigma$-coinvariant quotient torus of ${\rm T}\subseteq G$) by the definition of
the Kottwitz-Shelstad norm. See \cite{KS}, chapter 3 and \cite{Wes}, chapter 5. The relative
factor $vol_{db_\infty'}((G_{\gamma_0,\sigma}^I)'(\R))/
vol_{db_\infty'}((G_1)_{\gamma_1,id}^I)'(\R))$ turns out to be the degree of the isogeny
$T^{\theta^*} \to T_{\theta^*}$ independently from $I$. Notice, up to conjugacy over $\R$
(hence up to isomorphism over $\R$) the tori $T$ respectively $T_1$ only depend on $I$ and not
on $\gamma_0,\gamma_1$.

\bigskip\noindent
{\it Example}. For $G=(PGl( 2n+1),\sigma)$ and $G_1=(Sp(2n),id)$ as in \cite{Wes} (5.18) the
quotient $\vert \alpha_\infty(\gamma_0,1)\vert/\vert \alpha_\infty(\gamma_1,1)\vert$ is equal
to the relative measure factor $2^n$. In fact $K_\infty^{I,m}=K_\infty^{I,+}$ both for
$(G,\sigma)$ and $(G_1,id)$, since $O(2n+1,\R)/\{\pm id\} \cong SO(2n+1,\R)$ and since $U(n)$
is connected. If it holds for $I=\Delta$, then for all $I$. In all cases $H^1(\R,T)\cong
H^1(\R,T_1)$. For $M_I\cong PGl(2m+1)\times \prod_i Gl(r_i) \subseteq G= PGl(2n+1)$ and the
corresponding $(M_1)_I \cong Sp(2m)\times \prod_i Gl(r_i) \subseteq G_1=Sp(2n)$ only the cases
$r_i=0,1,2$ give nonvanishing contributions to the trace formula, using that $T$ and $T_1$ are
anisotropic modulo the center of $M_I$ resp. $(M_{1})_I$. The corresponding decomposition
$T=T_m \times \prod_i T_{r_i}$ and $T_1=T_{1,m}\times \prod_i T_{1,r_i}$ easily implies
$H^1(\R,T_m)\cong H^1(\R,T_{1,m})$ (being $\R$-anisotropic of the same rank) and
$H^1(\R,T_{r_i}) \cong H^1(\R,T_{1,r_i})$ by direct inspection, again using that $T$ and $T_1$
are anisotropic modulo the center of $M_I$ resp. $(M_{1})_I$.

\bigskip\noindent
We have shown

\bigskip\noindent
{\bf Lemma}. {\it For matching elements strongly $\sigma$-regular elements $\gamma_0$ and
matching $(\gamma_0,\gamma_1)$ the quotient
$\vert\alpha_\infty(\gamma_0)\vert/\vert\alpha_\infty(\gamma_1)\vert =c_I$ only depends on $I$
and $(G,\sigma,G_1)$, but not on the stable conjugacy class of $\gamma\in P_I(\Q)$. For
$I=\Delta$ this quotient defines a universal constant $c=c(G,\sigma,G_1)\neq 0$.}

\bigskip\noindent
In the other relevant case $(G,\sigma)=(GSO(3,3),\sigma)$ and $G_1=GSp(4)$ the constants $c_I$
depend on $I$.

\bigskip\noindent
\goodbreak\noindent
{\bf Theorem 4}. {\it For $G=GSO(3,3)$ resp. $G=PGl(5)$ and $\sigma$ as above, so that $G_1$ is
$GSp(4)$ respectively $G_1$ is $Sp(4)$, the following holds \smallskip \begin{enumerate}
\item{} Suppose $\Pi_1$ is a cuspidal automorphic
representation of $GSp(4,\A)$ (respectively an irreducible component of its restriction to
$Sp(4,\A)$), which is not CAP (nor a weak endoscopic lift). Suppose $(\Pi_1)_{fin}$ contributes
to $H^{\bullet}(M_{G_1},{\cal V}_1)$. Suppose $S$ is a sufficiently large finite set of places
for which $\Pi_1^S$ defined by $\Pi_1=\Pi_{1,S}\Pi_1^S$ is nonarchimedean and unramified. Let
$\Psi$ be the virtual representation of $G_1(\A_{fin})=G_1(\A_S)\times G_1(\A^S)$ on the
generalized (cuspidal) eigenspaces contained in $H^{\bullet}(M_{G_1},{\cal V}_1)$, on which
$G_1(\A^S)$ acts by $\Pi_1^S$. Then there exists a pair of globally and strongly matching
functions $f,f_1$ with strongly regular support, so that $Trace((\Pi_1)_{fin}(f_1)) \neq 0 $
and
$$Trace(\Psi (f_1)) \neq 0 \ .$$
%$$ and
%$$Trace((\Pi'_1)_f(f_1)) = 0 $$
%holds for all irreducible representations $\Pi'_1\not\cong \Pi_1$ of $G_1(\A_{fin})$ which
%contribute to $H^{\bullet}(M_{G_1},{\cal V}_1)$.
\item{} The fundamental lemma holds: In the situation
of 1) there exists a finite set of places $S=S(\Pi_1,G_1,G)$ such that for $v\not\in S$ the
representation $\Pi_{1,v}$ is unramified, the spherical Hecke algebra
$C_c^{\infty}(G(\Q_v)/\!/K_{v})$ is defined, so that for all $f_v\in C_c^{\infty}(G_v/\!/K_v)$
the functions $f_v$ and $b_{\xi}(f_v)$ are locally matching functions for the endoscopic datum
$(G_1,{}^{L}G_1,\xi,s)$ for $(G,\sigma,1)$ in the sense of \cite{KS}, (5.5.1) up to an
appropriate $z$-extension. So in particular the following generalized Shintani identities hold
for $v\notin S$
$$ Trace(b_{\xi}(f_v);\Pi_{1,v}) \ = \ Trace(f_v\cdot \sigma;
r_{\xi}(\Pi_{1,v})) \ .$$
\item{} For globally and strongly matching functions $(f,f_1)$ there exists a universal constant
$c=c(G,G_1)\neq 0$ such that the trace identity
$$ T(f,\sigma,G,\chi) = c\cdot T(f_1,id,G_1,\chi_1) \ $$
holds. \end{enumerate}}

\bigskip\noindent
{\it Proof of theorem 4}.  Assertion (4.3) has already been explained.

\bigskip\noindent
{\it Assertion (4.2)}. The comparison (4.3) of trace formulas would be pointless unless there
exist matching functions at least of the form (4.2). Assertion (4.2) is a statement, which is
enough to prove for almost all places $v\notin S$ for a suitable large finite set $S$. For an
arbitrary spherical Hecke operator at $v\notin S$ it can be reduced to the case, where
$f_v=1_{K_v}$ is the unit element of the spherical Hecke algebra. In the untwisted case this
was shown by Hales \cite{Ha} in full generality. In \cite{W3} this is done for a large class of
twisted cases, including those considered above, by extending a method of Clozel and Labesse.
We remark, that in the situation of theorem 4, the arguments of \cite{W3} can be simplified by
the use of the topological trace identity (4.3) explained above. Concerning unit elements: The
case of the unit element $f_v=1_{K_v}$ was established by Flicker \cite{Fl} for the first kind
of trace comparison involving $GSO(3,3)$ or more precisely $Gl(4)\times Gl(1)$, and was deduced
from that result for the second kind of trace comparison involving $PGl(5)$ in \cite{BWW}
theorem 7.9 and corollary 7.10. In both cases this is obtained for the unit elements of the
spherical Hecke algebras for large enough primes and for sufficiently regular
$(\gamma_0,\gamma_1)$. That these regularity assumptions do no harm is shown below.

\bigskip\noindent
{\it Assertion (4.1)}.  We claim that the assertions (4.2) and (4.3) of theorem 4 imply
assertion (4.1). Recall that $(\Pi_1)_{fin}$ and the coefficient system ${\cal V}_\chi$ are
fixed. Since $(\Pi_1)_{fin}$ is cuspidal but not CAP, it only contributes to cohomology in
degree 3, if it contributes nontrivially to the Euler characteristic of ${\cal V}_\chi$ (in our
case this amounts to the assumption that the archimedean component belongs to the discrete
series). All constituents of $\Psi$, being weakly equivalent to $\Pi_1$ and isomorphic to
$\Pi_1$ outside $S$, are not CAP. Hence the same applies for them.

\bigskip\noindent
To construct $f_1$ such that $T(f_1,id,G_1,\chi)\neq 0$ and $Trace(\Psi(f_1))\neq 0$ the
easiest candidate to come into mind is the following: Let $N$ be some principal congruence
level, i.e. assume $(\Pi_1)_{fin}^{K(N)}\neq 0$ where $K(N)\subseteq GSp(4,\Z_{fin})$ is
defined by the congruence condition $k\equiv id$ mod $N$. Choose a sufficiently large finite
set of places $S$ containing the finite set of divisors of $N$. Then put $f_1=\prod_{v\neq
\infty} f_{1,v}$, where $f_{1,v}$ is chosen to be the function which is zero for all $g\neq
k\cdot z$ for $k$ is in the principal congruence subgroup of level $N$ and chosen equal to
$f_{1,v}=\omega(z)^{-1}$ else. Then $f_{1,v}$ is the unit element of the Hecke algebra at
almost all places not in $S$ and the required condition $Trace((\Psi)_{fin}(f_1)) \geq
Trace((\Pi_1)_{fin}(f_1)) > 0$ holds by definition. To start to prove (4.1) we first have to
find a matching function $f$ on $G(\A_{fin})$. This involves the fundamental lemma, which is
known for elements $(\gamma_0,\gamma_1)$ sufficiently regular. Therefore we modify our first
naive choice $f_1$ slightly. For this we choose an auxiliary place $w\notin S$, where $\Pi_1$
is unramified. At this auxiliary place $w$ our previous function $f_{1,w}$ will be replaced by
an elementary function in the sense of \cite{La}, \cite{W3}. We will see, that this allows us
to restrict ourselves to consider matching conditions at sufficiently regular elements
$(\gamma_0,\gamma_1)$. Once these data in $S$ and $w$ are fixed, we  finally allow for some
additional modification at other unramified place $w'\neq w$ in order to construct a good
$\Pi_1^S$-projector in the sense of \cite{Lec}.

\bigskip\noindent
%It is therefore enough to find matching functions $f_v$ for the nonarchimedean places $v\in S$,
%since then $f=\prod_{v\neq \infty} f_v$ defines a pair $(f_1,f)$ of matching functions with the
%required property (4.1). In fact, we allow a minor modification at one of the unramified
%places. We assume $S$ to be large enough, so that it contains a nonarchimedean place $w$, where
%$\Pi_1$ is unramified.
At the auxiliary place $w$ we choose $f_{1,w}$ to be an elementary function as in \cite{La} \S
2 attached to some sufficiently regular element $t_0$ in the split diagonal torus, to be
specified later. Then there exists a matching $\sigma$-twisted elementary function on $G(\Q_w)$
\cite{W3}. Notice, the stable orbital integrals $SO^{G_1}_t(f_{1,w})$ vanish, unless $t$ is
conjugate in $G_1(F)$ to some element in the torus of diagonal matrices, which differs from
$t_0$ by a unimodular diagonal matrix and a central factor. In other words, the stable orbital
integral of $f_{1,w}$ then locally has regular semisimple support for suitable sufficiently
regular element $t_0$. For the global topological trace formulas of $G_1$ and $G$ this has the
effect, that only regular semisimple global elements give a nonzero contribution on the
geometric side. Hence for the construction a global matching functions $(f_1,f)$ it suffices to
show local matching for $(f_v,f_{1,v})$ at $v\neq w$ by considering stable orbital integrals at
regular elements only! Since in our situation the Kottwitz-Shelstad transfer factors are
identically one for all regular elements, this means we have to show local matching of stable
orbital integrals at all semisimple regular elements $\gamma_1$ at the places $v\neq w$. In
addition we have to find a function $f_w$ matching with $f_{1,w}$ at the place $w$ supported in
strongly $\sigma$-regular elements $\gamma_0$ of $G(\Q_w)$, and we have to guarantee
$Trace(\Psi_{1,w})(f_{1,w})\neq 0$ by a suitable choice of $t_0$. By taking finite linear
combination of elementary functions of this type one can also achieve the vanishing of all
representation $Trace(\pi_{1,w})(f_{1,w})$ for the finitely many representations $\pi=\Pi_i$
weakly equivalent to $\Pi_1$ with level $N\cdot p_w$, which contribute to cohomology of the
fixed coefficient system ${\cal V}_\chi$. Hence the trace coincides with the trace of our first
naive choice of $f_1$. In fact, as already mentioned, one can then also find a corresponding
finite linear combination $f_{w}$ of $\sigma$-twisted elementary functions matching with
$f_{1,w}$, and for this choice of $f_{1,w}$ then $f_w$ has support in strongly $\sigma$-regular
elements (see \cite{W3}).

\bigskip\noindent
This modification using $t_0$ at the place $w$, which we further discuss below, is useful also
for other purposes. At the moment, taking the situation at the place $w$ for granted, let us
look at the other places first. For the nonarchimedean places $v\neq w$ not in $S$ we now can
use assertion (4.2) to obtain the matching functions $f_v$ to be unit elements at all these
places. Recall, the fundamental lemma for the unit elements at sufficiently regular suffices
for this. Concerning the place $v\in S$, the existence of local matching function  can be
reduced to existence of matching germs. This can be achieved by the argument of \cite{LS}, \S
2.2 or the similar argument of \cite{V}. So one is reduced to the matching of germs of stable
orbital integrals.  Since we deal with a purely local question now, of course one has to allow
singular elements now. Fortunately for the pairs $(G,G_1)$ under consideration the matching of
germs of stable orbital integrals has been proven by Hales \cite{Ha2}. This completes our
construction of the functions $f_v$ at the places $v\neq w$ defining the globally matching pair
$(f_1,f)$. So let us come back to the auxiliary unramified place $w$.

\bigskip\noindent
Concerning the choice of $t_0$: The trace $Trace(f_1,\Pi_1)=0$ of any irreducible cuspidal
automorphic representation $\Pi_1$ vanishes unless $\Pi_1$ admits a nontrivial fixed vector for
the group $K(N,w)=K_w\prod_{v\neq \infty,w}K(N)_v$, where $K_w$ is chosen to be a Iwahori
subgroup. Hence only finitely many irreducible automorphic representations with fixed
archimedean component and $Trace(f_1,\Pi_1)\neq 0$ exist.  $Trace(\Pi_{fin}(f_1))$ considered
in assertion (4.1) is a linear combination of the type $\sum_{i} m(\Pi_i) Trace(f_{1,w},\Pi_i)$
as a function of $t_0$, involving a finite number of local irreducible admissible
representations $\pi=\Pi_i$ of $G_1(\Q_w)$ of Iwahori level, with certain multiplicities
$m(\Pi_i)
> 0$. Since the unramified representation $\Pi_{1,w}$ is one of these, this sum is
nonempty. By the argument of \cite{La} (section 5)  $t_0$, which defines the elementary
function $f_{1,w}$, can be chosen to be strongly regular so that this sum is nonzero, provided
$G_1$ is of adjoint type. This adjointness assumption is true for $G_1=GSp(4)$, but not for
$G_1=Sp(4)$. However, since we consider irreducible components $\Pi_1$ of the restriction of
the global representation in the case of the group $Sp(4)$, this kind of argument carries over
also for the case $G_1=Sp(4)$. Of course this remark proves the required nonvanishing
$Trace(\Pi(f_{1})) \neq 0$, since it allows reconstruct the spherical trace of the spherical
function of $\Pi_{1,w}$ at the place $w$ by a linear combination $f_{1,w}$ of elementary
functions attached to finitely many strongly regular elements $t_0$.

\bigskip\noindent
We finally make some further adjustment, but only using spherical functions. This is necessary
to obtain strongly matching functions, and this is only necessary in the case
$(G,G_1)=(GSO(3,3),GSp(4))$. For this we consider the unramified representation $\Pi'^S =
r_\xi(\Pi^S)$. Using the method of CAP-localization \cite{Lec} (with levels and coefficient
system fixed) it is clear that we can modify our $f$ at unramified places not in $S$ using a
good spherical $\Pi'^S$-projector $\tilde f^S$. Similarly modify $f_1$ by the corresponding
good spherical $\Pi_1^S$-projector $\tilde f_1^S=b_\xi(\tilde f^S)$. Then all contributions
$T_I(f_1,id,G_1,\chi)$ in the trace formula for $P_I\neq G_1$ in the trace formula vanish. This
is possible, since we assumed $\Pi^S$ not to be a CAP-representation. Then, since $f$ is
globally matching, we get $T_I(f,G,\sigma,\chi) = c_I \cdot T_I(f_1,id,G_1,\chi)$ for constants
$c_I$ by the last lemma. Hence also $T_I(f,G,\sigma,\chi)=0$ for $P_I\neq G$. Therefore the
modified pair $(f\tilde f^S,f_1\tilde f_1^S)$ is not only globally matching, but also a
strongly matching pair of functions with the properties required by assertion 4.1. This
completes the proof of theorem 4.

\bigskip\noindent
{\it Applying theorem 4}. For an irreducible cuspidal representation $\Pi$, which is not CAP,
suppose $\Pi_{\infty}$ belongs to the discrete series. For $G_1(\A)=GSp(4,\A)$ the
representation $\Pi_{fin}$ contributes to $H^3(M_1,{\cal V}_{\chi_1})$ for some character
$\chi_1$ ( \cite{Ast}, hypothesis $A$ and $B$) and  we can therefore assume that
$Trace(\Pi_{fin}(f_1);H^{\bullet}(M_1,{\cal V}_{\chi_1}))$ does not vanish for all $f_1$. Now
choose $f_1$ and a globally and strongly matching function $f$ on the group $G(\A_{fin})$ for
$G=GSO(3,3)$ (or its z-extension $Gl(4)\times Gl(1)$) as in theorem (4.1). Then by theorem
(4.2) and (4.3) there must exist a $\sigma$-stable irreducible $G(\A_{fin})$-constituent
$(\Pi',\omega)$ of $H^{\bullet}(M_G,{\cal V}_{\chi})$, for which $Trace(\Pi',\omega)(f\cdot
\sigma)\neq 0$ holds. By a theorem of Franke \cite{Fr} this representation is automorphic. Any
modification of the pair $(f_1,f)$ to $(f_{1,v} \prod_{w\neq v} f_{1,w}, f_v\prod_{w\neq v}
f_{w} )$ at a place $v\not\in S$ for $f_{1,v}=b_{\xi}(f_v)$ and spherical Hecke operators $f_v$
again gives a pair of globally strongly matching functions. Then $(4.2)$ and $(4.3)$ and the
linear independence of characters of the group $G(\A_{fin})\rtimes <\sigma>$ imply, that the
representation identity $(\Pi'_v,\omega_v)=r_{\xi}(\Pi_v)$ holds for all $v\not\in S$. Since
the partial $L$-series outside $S$ can be easily computed from the Satake parameters, we have
constructed a weak lift: The $L$-identity
$$ L^S(\Pi'\otimes\chi,s)=L^S(\Pi\otimes\chi,s) $$ holds for all idele
class characters $\chi$ for a sufficiently large finite set $S$. This is just a transcription
of the fundamental lemma, i.e. follows from  part (4.2) and (4.3) of theorem 4, and defines the
weak lift $\Pi'$. The weak lift $\Pi'$, now constructed, is uniquely determined by $\Pi$, since
by the strong multiplicity 1 theorem for $Gl(n)$ it is determined by the partial $L$-series
$L^S(\Pi'\otimes\chi,s)$ (at the places $v\not\in S$), once it exists. This shows, that our
second temporary assumption is satisfied. The required identity for the central characters
$\omega_{\Pi}^2 = \omega_{\Pi'}$ holds, since it locally holds  at all nonarchimedean places
outside $S$ by the $L$-identity above.

\bigskip\noindent
The case $G_1=Sp(4),G=PGl(5)$ is completely analogous and the corresponding considerations now
also allows us to get rid of our first temporary assumption.

\bigskip\noindent
{\it Concerning the archimedean place}. The topological trace formulas compute  the traces of
the Hecke correspondences on the virtual cohomology of a given coefficient system in terms of
orbital integrals. Underlying the trace comparison (4.3) is a corresponding fixed lift of the
coefficient system ${\cal V}_{\mu_1}$. The coefficient system being fixed, only finitely many
archimedean representations contribute to the topological trace formula for this coefficient
system. Also notice, the representation $\Pi'_{\infty}$ has to be a representation of
$Gl(4,\R)$, which has nontrivial cohomology for the given lift of the coefficient system. In
fact this determines $\Pi'_{\infty}$ in terms of ${\cal V}_{\mu}$ and the cohomology degree.
Furthermore in our relevant case, i.e. for the lift to $Gl(4)\times Gl(1)$, well known
vanishing theorems for Lie algebra cohomology and duality completely determine $\Pi'_{\infty}$
in terms of the coefficient system ${\cal V}_{\mu}$. We leave this as an exercise. This
uniquely determines $\Pi'_\infty$ in terms of $\Pi_\infty$. This yields matching archimedean
representations, which have nontrivial cohomology in degree 3 for the group $H$ with
$\sigma$-equivariant representations, that have nontrivial cohomology for $G$.

\bigskip\noindent
Consider the theta lift from $GO(3,3)$ to $GSp(4)$, which maps the class of $(\Pi',\omega)$ to
the class of $\Pi_{gen}$. At the archimedean place we claim, that for this lift
$(\Pi'_\infty,\omega_\infty)$ locally (!) uniquely determines $\Pi_{gen,\infty}$. Since theta
lifts behave well with respect to central characters, it is enough to prove this for the dual
pair $(O(3,3),Sp(4))$ by our earlier remarks on $\Pi'_\infty$. In fact by Mackey's theory the
restriction of $\Pi'_\infty$ to $Sl(4,\R)\cdot \R^*$ decomposes into two nonisomorphic
irreducible constituents, from which (from each of them) in turn $\Pi'_\infty$ is obtained  by
induction. The same holds for the restriction of discrete series representations of $GSp(4,\R)$
to $Sp(4,\R)\cdot \R^*$. It is therefore enough to determine the restriction of
$\Pi_{gen,\infty}$ to $Sp(4,\R)\cdot \R^*$. Once we have shown, that this restriction contains
an irreducible constituent in the discrete series, we have therefore as a consequence, that the
representation $\Pi_{gen,\infty}$ of $GSp(4,\R)$ is uniquely determined by $\Pi'_\infty$ and
that $\Pi_{gen,\infty}$ belongs to the discrete series.  We now still have to understand, why
$\Pi_{gen,\infty}$ should be in the archimedean $L$-packet of $\Pi_\infty$.

\bigskip\noindent
This being said, we first replace $Gl(4,\R)$ by $Sl(4,\R)$ or the quotient
$SO(3,3)(\R)=Sl(4,\R)/\Z_2$ (recall the central character was $\omega^2$) and replace
$\Pi'_\infty$ by a suitable irreducible constituent of the restriction. For simplicity of
notation still denote it $\Pi'_\infty$. If the representation on $O(3,3)(\R)$ induced from the
irreducible representation $\Pi'_\infty$ of $SO(3,3)(\R)$ would be irreducible, we can
immediately apply \cite{H} to conclude, that its theta lift on $Sp(4,\R)$ is uniquely
determined by $\Pi'_\infty$. Otherwise there exist two different extensions of $\Pi'_\infty$ to
$O(3,3)(\R)$, by simplicity denoted $\Pi'_\infty$ and $\Pi'_\infty\otimes \varepsilon$, where
$\varepsilon$ is the quadratic character of $O(3,3)(\R)$ defined by the quotient
$O(3,3)(\R)/SO(3,3)(\R)$. However if this happens, we now claim only one of the two
possibilities contributes to the theta correspondence, so that again we can  apply \cite{H}.
The claim made follows from an archimedean version of part c) of the Proposition of \cite{V},
p.483. For the convenience of the reader we prove this in the archimedean exercise below.
\goodbreak

\bigskip\noindent
%By looking at weak endoscopic lifts $\Pi$, we see that
%$\Pi'_{\infty}$ is induced from a parabolic subgroup whose Levi
%subgroup contains $GO(2,2)(\R)$. Hence the study of the theta
%lift at the archimedean place can be reduced to a study of the
%easier lift from $GO(2,2)$ to $GSp(4)$. We skip the details and
%refer to [Wa], or [W2], p.5ff.
So passing from $\Pi_{\infty}$ to $(\Pi'_{\infty},\omega_{\infty})$ and then back to
$\Pi_{gen,\infty}$ turns out to be a well defined local assignment at the archimedean place; to
be accurate, in the first instance only up to twist by the sign-character. That it is well
defined then follows a posteriori, once we know the assigned image is contained in the discrete
series. For this see the archimedean remarks made following our second temporary assumption. To
compute the local assignment - first only up to a possible character twist - it is enough to do
this for a single suitably chosen global automorphic cuspidal representation $\Pi$. We have to
show, that passing forth back with these two lifts locally at the archimedean place, we do not
leave the local archimedean $L$-packet. Since for every archimedean discrete series
representation $\Pi_{\infty}$ of weight $(k_1,k_2)$ there exists a global weak endoscopic lift
$\Pi$ with this given archimedean component $\Pi_\infty$ (see \cite{Ast}), it is now enough to
do this calculation globally for this global weak endoscopic lift $\Pi$. For global weak
endoscopic lifts there exists a unique globally generic cuspidal representation $\Pi_{gen}$ in
the weak equivalence class of $\Pi$ by \cite{Ast}, hypothesis A. So $\Pi_{gen}$ is uniquely
determined by the description of this generic component given in \cite{Ast}, hypothesis A.
Therefore looking at the archimedean place the archimedean component $\Pi_{gen,\infty}$, which
is uniquely determined by $\Pi_\infty$ as we already have seen, has to be this unique generic
nonholomorphic member of the discrete series $L$-packet of $GSp(4,\R)$ of weight $(k_1,k_2)$.
Since it depends only on the archimedean local representation, and since this holds true in the
special case, this holds in general. This proves theorem 1 modulo the following

\bigskip\noindent
{\it Archimedean exercise}. We now prove the archimedean analog of proposition \cite{V}, p.483,
which we used above. See also \cite{Wei}. Let $V$ be a nondegenerate real quadratic space of
dimension $m=p+q$, let $G=O=O(p,q)$ be its orthogonal group with maximal compact subgroup
$K=O(p)\times O(q)$. Let $O^{\vee{}}$ be the set of equivalence classes of irreducible
$(Lie(O),K)$-modules. The metaplectic cover $Mp(2N)$ of $Sp(2N)$ for $N = m\cdot n$ naturally
acts on the Fock space. Its associated Harish-Chandra module $P_V(n) \cong
Sym^{\bullet}(V\otimes_{\R} \C^n)$ defines the oscillator representation $\omega_{Fock}$. For
the pair $(G,G')=(O,Mp(2n,\R))$ the restriction of $\omega_{Fock}$ to $G\times G' \subseteq
Mp(2N,\R)$ induces the theta correspondence. The action of $G$ by $\omega_{Fock}$ does not
coincide with the natural action on $Sym^{\bullet}(V\otimes_{\R} \C^n)$, except when $V$ is
anisotropic. Nevertheless the action of $K$ does. Let $R(n)\subseteq O^{\vee{}}$ denote the
classes of irreducible $(Lie(O),K)$-module quotients of $P_V(n)$. For $\pi\in O^{\vee{}}$ let
$n(\pi)$ be the smallest integer $n$, if it exists, such that $\pi\in R(n)$. Since
$P_V(n+n')=P_V(n)\otimes_{\C} P_V(n')$  - considered as modules of $(Lie(O),K)$ - for any
$\pi\in R(n)$ and $\pi'\in R(n')$ the irreducible quotients of $\pi\otimes_{\C}\pi'$ contribute
to $R(n+n')$. For $\pi\in O^{\vee}$ let $\pi^{\vee}$ denote the contragredient representation.
Since the quadratic character $\varepsilon$ is an irreducible quotient of the representation
$\pi\otimes (\pi^{\vee}\otimes\varepsilon)$, we obtain $n(\pi) + n(\pi^{\vee}\otimes
\varepsilon) \geq n(\varepsilon)$. We claim $n(\varepsilon)\geq dim_{\R}(V)$, which implies as
desired
$$ n(\pi) + n(\pi^{\vee}\otimes \varepsilon) \ \geq \ dim_{\R}(V) \ .$$

\bigskip\noindent
For the proof of $n(\varepsilon)\geq m$ put $n=n(\varepsilon)$. The restriction of
$\varepsilon$ to $K=O(p)\times O(q)$ is $\sigma=\varepsilon\boxtimes\varepsilon$.
$M'=Mp(2n,\R)\times Mp(2n,\R)$ covers the centralizer of $K$ in $Sp(2N)$ and $M'_0=\tilde
U(n,\C)\times \tilde U(n,\C)$ its maximal compact subgroup. By \cite{KV} and \cite{H}, lemma
3.3 a unique irreducible representation $\tau'$ of $M'_0$ is attached to $\sigma$, which is the
external tensor product of the irreducible highest weight representations with highest weights
$({p\over 2}+1,..,{p\over 2}+1,{p\over 2},...,{p\over 2})$ respectively $({q\over
2}+1,..,{q\over 2}+1,{q\over 2},...,{q\over 2})$ of $\tilde U(n,\C)$, which is a twofold
covering group of $U(n,\C)$. Here ${p\over 2}+1$ occurs $p$ times respectively ${q\over 2}+1$
occurs $q$ times. So in particular $n\geq max(p,q)$. The two components are represented by
pluriharmonic polynomials in the Fock space of degree $p$ respectively $q$. As a consequence,
$\tau'$ is realized in $Sym^d(V\otimes_{\R} \C^n)$ for $d=p+q=m$. In other words $deg(\tau')=m$
for the degree $deg$ in the sense of \cite{H}. The restriction of $\tau'$ to the maximal
compact subgroup $K' = \tilde U(N,\C) \cap G'$ is the tensor representation $det^{{p-q \over
2}}\otimes \bigwedge^p(\C^n) \otimes \bigwedge^q(\C^n)^{\vee}$ of $\tilde U(n,\C)$, since the
negative definite part gives an antiholomorphic embedding (\cite{H}, p.541). The highest
weights of all irreducible constituents of the representation $\bigwedge^p(\C^n) \otimes
\bigwedge^q(\C^n)^{\vee}$ of $U(n,\C)$ are of the form $(1,..,1,0,..,0,-1,...,-1)$ with $i\leq
p$ digits $1$, $j\leq q$ digits $-1$ and $n-i-j$ digits $0$. According to \cite{H}, lemma 3.3
there must be a unique irreducible representation $\sigma'$ of $\tilde U(n,\C)$ in the
restriction of $\tau'$ such that $deg(\sigma')=deg(\tau')=m$. The representation of $M_0'\cong
\tilde U(n,\C)\times \tilde U(n,\C)$ on the polynomials of degree $k$ in the Fock space
$P_V(n)$ is isomorphic to $\bigoplus_{a+b=k} det^{{p\over 2}}\cdot Sym^a(\R^p\otimes_{\R} \C^n)
\boxtimes det^{{q\over 2}}\cdot Sym^b(\R^q\otimes_{\R} \C^n) $. Its restriction to $K'\cong
\tilde U(n,\C)$ in $G'$ therefore is isomorphic to the representation $det^{{p-q\over 2}}
\bigoplus_{a+b=k} Sym^a(\R^p\otimes_{\R} \C^n) \otimes Sym^b(\R^q\otimes_{\R} (\C^n)^{\vee})$
induced by the natural action of $U(n,\C)$ on $\C^n$. Hence for $i\leq p$ and $j\leq q$ (the
number of digits $1$ resp. $-1$ of the highest weight) we get $deg(det^{{p-q\over 2}}\otimes
(1,..,1,0,0,..,-1,..,-1)) = i+j $ by considering the representation generated by the product of
some $i\times i$-minor in $Sym^i(\R^i\otimes_{\R} \C^n)$ and some $j\times j$-minor in
$Sym^i(\R^j\otimes_{\R} \C^n)$. Thus $deg(\sigma')=m$ implies, that there exists $i\leq p$ and
$j\leq q$ such that $deg(\sigma')=i+j=m$. Therefore $i=p$ and $j=q$, hence $m=i+j\leq n$. So we
have shown $n(\varepsilon)\geq m$. In fact it is not hard to see $n(\varepsilon)=m$.
%For this we use the two equivalent descriptions of the theta
%correspondence by [H], theorem 1 in terms of the smooth
%representation $\omega^{\infty}$ and equivalently by the
%algebraic formulation [H], theorem 2.1 in terms of the Fock model
%$\omega_{Fock}$. From the formulation in terms of the smooth
%model $\varepsilon\in R(k)$ implies existence of a nontrivial
%tempered distribution $T$ on $V^k$ with $T^g =
%\varepsilon(g)\cdot T$ for all $g\in O$. Note that the space of
%smooth vectors of the unitary oscillator representation $\omega$
%of $Sp(k\cdot dim_{\R}(V))$ considered as a topological
%vectorspace is the Schwartz space ${\cal S}(V^k)$ with the
%Schwartz topology ([BW], lemma 1.11 on page 153). By the argument
%in the appendix of [R], which also applies in our present
%archimedean setting, such a distribution $T$ vanishes unless
%$k\geq dim_{\R}(V)$.

\bigskip\noindent
{\it The special case considered}. In the situation of theorem 1 and its proof this gives
$$n(\pi) + n(\pi\otimes\varepsilon) \geq 6\ .$$ The underlying representation $\Pi'_\infty$ of
$GO(3,3)(\R)$ decomposes into two nonisomorphic representations $\pi_1,\pi_2$ of $O(3,3)(\R)$
with $n(\pi_1)=n(\pi_2)=2$. Therefore $n(\pi_1\otimes\varepsilon) \geq 4$ and
$n(\pi_2\otimes\varepsilon) \geq 4$ by the above inequality, since in our case $\pi^\vee \cong
\pi$ follows from $(\Pi'_{\infty})^{\vee} \cong \Pi'_{\infty}\otimes\omega_{\infty}^{-1}$
(which was a consequence of the second temporary assumption).

\bigskip\noindent
This finally completes the proof of theorem 1.

\bigskip\noindent
\goodbreak
\bigskip\noindent
{\it Remark on orthogonal representations}. For a weak endoscopic lift $\Pi$ the statement of
theorem 1 is known by \cite{Ast}, hypothesis A. So there was no need not go to the trace
formula arguments required otherwise. Nevertheless the above arguments applied  for a weak
endoscopic lift nevertheless yield something interesting. Let $\Pi$ be cuspidal irreducible
representation $\Pi$ of $GSp(4,\A)$, which is a weak endoscopic lift attached to a pair of
cuspidal representations $(\sigma_1,\sigma_2)$ of $Gl(2,\A)$ with central character $\omega$,
but which is not CAP. Then we can still construct the irreducible automorphic representation
$\tilde\Pi'$ of $Gl(5,\A)$ from $\Pi$ as above. However the representation $\tilde\Pi'$ will
not be cuspidal any more, since its $L$-series $L^S(\tilde\Pi',s)= L^S(\sigma_1 \times \sigma_2
\otimes\omega^{-1},s)\zeta^S(s)$ has a pole at $s=1$. Hence the automorphic representation
$\tilde\Pi'$ is Eisenstein, in fact induced from an automorphic irreducible representation
$(\pi_4,\pi_1)$ of the Levi subgroup $Gl(4,\A)\times Gl(1,\A)$. This gives rise to an
irreducible automorphic representation $\pi_4$ of $Gl(4,\A)$ attached to $(\sigma_1,\sigma_2)$,
which after a character twist by $\omega$ will be denoted
$$ \sigma_1\times \sigma_2 \ .$$ The automorphic representation
$\sigma_1\times\sigma_2$ is uniquely determined by its partial $L$-series
$L^S(\sigma_1\times\sigma_2,s)$, which coincides with the partial $L$-series attached to the
orthogonal four dimensional $\overline E_\lambda$-valued $\lambda$-adic Galois representation
$\rho_{\sigma_1,\lambda}\otimes_{\overline E_{\lambda}} \rho_{\sigma_2,\lambda}$, the tensor
product of the two dimensional Galois representations $\rho_{\sigma_i,\lambda}$ attached to
$\sigma_1$ and $\sigma_2$. This four dimensional tensor product is an orthogonal Galois
representation. It should not be confused with the symplectic four dimensional Galois
representation obtained in \cite{Ast}, theorem I, which for weak endoscopic lifts is the direct
sum of the Galois representations $\rho_{\sigma_i,\lambda}$. D.Ramakrishnan obtained a
completely different and more general construction of the automorphic representation
$\sigma_1\times \sigma_2$ using converse theorems.

\goodbreak
\bigskip\noindent
\medskip\noindent

\bigskip\noindent
\bigskip\noindent\footnotesize
\centerline{\rm Rainer Weissauer} \centerline{ \rm Mathematisches Institut, Im Neuenheimer Feld
288} \centerline{\rm Universit\"at Heidelberg, 69120 Heidelberg} \centerline{\rm email:
weissauer@mathi.uni-heidelberg.de} \vfill\eject
\end{document}